\documentclass[11pt]{article}
\usepackage{amsmath}
\usepackage{amssymb}
\usepackage{array}
\newtheorem{theo}{\hspace*{\parindent}Theorem}

\newtheorem{lemma}{\hspace*{\parindent}Lemma}

\newtheorem{corolt}{\hspace*{\parindent}Corollary}[theo]
\def\sign{\mathrm{sign}}
\newcounter{theremark}
\newcommand{\rem}{\par\refstepcounter{theremark}\textbf{Remark \arabic{theremark}.} }
\title{An approximation for zero-balanced Appell function $F_1$ near $(1,1)$}
\author{D. Karp\footnote{Institute of Applied Mathematics, Vladivostok, Russia,
e-mail:\emph{dmkrp@yandex.ru}}}
\date{}
\begin{document}
\maketitle

\begin{center}
\parbox{12cm}{
\small\textbf{Abstract.} We suggest an approximation for the
zero-balanced Appell hypergeometric function $F_1$ near the
singular point $(1,1)$. Our approximation can be viewed as a
generalization of Ramanujan's approximation for zero-balanced
${_2F_1}$ and is expressed in terms of ${_3F_2}$. We find an error
bound and prove some basic properties of the suggested
approximation which reproduce the similar properties of the Appell
function. Our approximation reduces to the approximation of
Carlson-Gustafson when the Appell function reduces to the first
incomplete elliptic integral.}
\end{center}

\bigskip

\paragraph{1. Introduction.}  The generalized hypergeometric function
is defined by \cite[formula 4.1(1)]{Bat1}
\begin{equation}\label{eq:pFq}
{_{p}F_{q}}\left(\left.\!\!\begin{array}{l} a_1,\ldots, a_p
\\ b_1,\ldots, b_q
\end{array}\right|z\!\right)=\sum\limits_{k=0}^{\infty}\frac{(a_1)_k\ldots(a_p)_k}
{(b_1)_k\ldots(b_q)_k}\frac{z^k}{k!},
\end{equation}
where $(a)_0=1$, $(a)_k=a(a+1)\cdots(a+k-1)$, $k=1,2,\ldots$, is
shifted factorial.  This function is called zero-balanced if
$p=q+1$ and $\sum_{i=1}^{p}a_i=\sum_{i=1}^{q}b_i$.

Ramanujan (see \cite{Berndt0,Berndt1,Berndt2}) suggested the
following approximations for zero-balanced ${_2F_1}$  and
${_3F_2}$:
\begin{equation}\label{eq:2F1}
B(a,b){_2F_1}(a,b;a+b;x)=-\ln(1-x)+\gamma(a,b)+O((1-x)\ln(1-x)),~~x\to{1-},
\end{equation}
where
\begin{equation}\label{eq:constants}
\gamma(a,b)=2\psi(1)-\psi(a)-\psi(b),~~\psi(z)=\frac{\Gamma'(z)}{\Gamma(z)},
\end{equation}
 and
\[
\frac{\Gamma(a_1)\Gamma(a_2)\Gamma(a_3)}{\Gamma(b_1)\Gamma(b_2)}
{_{3}F_{2}}\left(\!\!\left.\begin{array}{l} a_1, a_2, a_3
\\ b_1, b_2
\end{array}\right|x\!\right)=-\ln(1-x)+L+O((1-x)\ln(1-x)),~~x\to{1-},
\]
where $\Re(a_3)>0$ and
\[
L=2\psi(1)-\psi(a_1)-\psi(a_2)+\sum\limits_{k=1}^{\infty}\frac{(b_2-a_3)_k(b_1-a_3)_k}{k(a_1)_k(a_2)_k}.
\]
These formulas have been generalized to ${_{q+1}F_q}$ by
N$\o$rlund \cite{Norlund}, Saigo and Srivastava in \cite{Saigo},
Marichev and Kalla in \cite{Marichev} and B\"{u}hring in
\cite{Buhring}, see details in the survey paper by B\"{u}hring and
Srivastava \cite{Srivastava}.

The Appell function $F_1$ generalizes ${_2F_1}$ to two variables
and  is defined by \cite{Bat1}:
\begin{equation}\label{eq:F1exp}
F_1(\alpha;\beta_1,\beta_2;\gamma;z_1,z_2)=\sum\limits_{k,n=0}^{\infty}
\frac{(\alpha)_{k+n}(\beta_1)_k(\beta_2)_n}{(\gamma)_{k+n}k!n!}z_1^kz_2^n,
\end{equation}
for $|z_1|<1$, $|z_2|<1$ and by analytic continuation for other
values of $z_1$, $z_2$. An asymptotic expansion for $F_1$ has been
studied by Ferreira and L\'{o}pez in \cite{Fer-Lopez} in the
neighborhood of infinity. This approximation can be converted into
an approximation around $(1,1)$ using the formula
\[
F_1(a;b,c;d;1-xz,1-yz)=x^{-b}y^{-c}F_1\left(d-a;b,c;d;1-\frac{1}{xz},1-\frac{1}{yz}\right).
\]
It has been noticed by B.C.\,Carlson in \cite{Carlson} that the
incomplete elliptic integral of the first kind is a particular
case of $F_1$:
\begin{equation}\label{eq:EllipticF-F1}
F(\lambda,k)={\lambda}F_1(1/2;1/2,1/2;3/2;\lambda^2,k^2\lambda^2).
\end{equation}
Carlson and Gustafson studied the asymptotic approximation for
$F(\lambda,k)$ in \cite{CG1}.  Their expansion can be shown to be
a particular case of the expansion for $F_1$ given later in
\cite{Fer-Lopez}. We will show below that both expansions (but not
the error bounds!) can be obtained by simple rearrangement of
(\ref{eq:F1exp}) and use of known transformation formulas for
$F_1$.  More precise approximations for $F(\lambda,k)$ which
cannot be reduced to expansions from \cite{Fer-Lopez} have been
given recently  by S.M.\,Sitnik and the author in
\cite{KarpSitnik}.

The purpose of this paper is to give an analogue of (\ref{eq:2F1})
for the "zero-balanced" Appell function $F_1$ with
$\gamma=\alpha+\beta_1+\beta_2$. Important properties of $F_1$ are
permutation symmetry
\begin{equation}\label{eq:F1sym}
F_1(\alpha;\beta_1,\beta_2;\gamma;z_1,z_2)=F_1(\alpha;\beta_2,\beta_1;\gamma;z_2,z_1),
\end{equation}
reduction formulas
\begin{equation}\label{eq:F1sidereduction}
F_1(\alpha;\beta_1,\beta_2;\gamma;z,1)={_2F_1}(\alpha,\beta_2;\gamma;1){_2F_1}(\alpha,\beta_1;\gamma-\beta_2;z),
\end{equation}
\begin{equation}\label{eq:F1digreduction}
F_1(\alpha;\beta_1,\beta_2;\gamma;z,z)={_2F_1}(\alpha,\beta_1+\beta_2;\gamma;z),
\end{equation}
and reduction formula (\ref{eq:EllipticF-F1}).  Our approximation
reproduces the permutation symmetry (\ref{eq:F1sym}), reduces to
Ramanujan approximation given in (\ref{eq:2F1}) in cases given by
(\ref{eq:F1sidereduction}) and (\ref{eq:F1digreduction}) and
reproduces Carlson-Gustafson approximation for the values of
parameters given in (\ref{eq:EllipticF-F1}).

Some new reduction formulas for $F_1$ have been discovered in
\cite{Ismail}.

\paragraph{2. Main results.} To save space let us introduce the notation
\begin{equation}\label{eq:f-defined}
f_{a,b_1,b_2}(x,y)=B(a,b_1+b_2)F_1(a;b_1,b_2;a+b_1+b_2;x,y).
\end{equation}
Our main approximation is given by
\begin{equation}\label{eq:g-defined}
g_{a,b_1,b_2}(x,y)=\ln\frac{1}{1-x} +\gamma(a,b_1+b_2)
+\frac{b_2(y-x)}{(b_1+b_2)(1-x)}
{_3F_2}\left(\!\!\begin{array}{l}1,1,b_2+1\\2,b_1+b_2+1\end{array}\!\vline~\frac{y-x}{1-x}\right),
\end{equation}
where $\gamma(a,b_1+b_2)$ is defined in (\ref{eq:constants}). The
following theorem confirms that $g_{a,b_1,b_2}$ is indeed a
correct analogue of the righthand side of (\ref{eq:2F1}).
\begin{theo}\label{th:fgapprox}
For $0\leq{x}<{1}$, $0\leq{y}<{1}$, $a,b_1,b_2>0$:
\begin{equation}\label{eq:g-approx}
f_{a,b_1,b_2}(x,y)=g_{a,b_1,b_2}(x,y)+R_{a,b_1,b_2}(x,y),
\end{equation}
with
\begin{equation}\label{eq:R1bound}
0<R_{a,b_1,b_2}(x,y)<r\left(1+a-a\ln(r)\right)=O(r\ln(r)),
\end{equation}
where in the last formula $x,y\to{1}$, $r=(1-x)b_1+(1-y)b_2\to{0}$
is the ''rhombic'' distance to $x=y=1$ which is asymptotically
equivalent to Euclidian distance.
\end{theo}
\begin{corolt} Formulas \emph{(\ref{eq:g-approx})} and \emph{(\ref{eq:R1bound})}  imply in particular the
inequality
\begin{equation}\label{eq:fgineq}
f_{a,b_1,b_2}(x,y)>g_{a,b_1,b_2}(x,y)
\end{equation}
for all $x,y\in(0,1)$.
\end{corolt}

\textbf{Proof of Theorem~\ref{th:fgapprox}.} A simple
rearrangement of (\ref{eq:F1exp}) gives
\begin{equation}\label{eq:F12F1}
F_1(\alpha;\beta_1,\beta_2;\gamma;z_1,z_2)
=\sum\limits_{k=0}^{\infty}\frac{(\alpha)_{k}(\beta_1)_k}{(\gamma)_{k}k!}{_2F_1}(\alpha+k,\beta_2;\gamma+k;z_2)z_1^k.
\end{equation}
Suppose $\gamma=\alpha+\beta_2$, then ${_2F_1}$ in
(\ref{eq:F12F1}) is zero-balanced and we can apply \cite[formula
2.10(12)]{Bat1}
\begin{multline}\label{eq:2F11mainlog}
\frac{\Gamma(\eta)\Gamma(\beta)}{\Gamma(\eta+\beta)}{_2F_1}(\eta,\beta;\eta+\beta;z)=
\sum\limits_{n=0}^{\infty}\frac{(\eta)_n(\beta)_n}{(n!)^2}
[-\log(1-z)+2\psi(n+1)-\psi(\eta+n)-\psi(\beta+n)](1-z)^n.
\end{multline}
It  gives
\begin{multline}\label{eq:F1logarithmic}
\frac{\Gamma(\alpha)\Gamma(\beta_2)}{\Gamma(\alpha+\beta_2)}F_1(\alpha;\beta_1,\beta_2;\alpha+\beta_2;z_1,z_2)=
\\
\sum\limits_{n,k=0}^{\infty}
\frac{(\alpha+k)_n(\beta_2)_n(\beta_1)_kz_1^k}{(n!)^2k!}
[-\ln(1-z_2)+2\psi(1+n)-\psi(\beta_2+n)-\psi(\alpha+k+n)](1-z_2)^n.
\end{multline}
Taking account of
\[
(\alpha)_{k+n}=(\alpha)_k(\alpha+k)_n=(\alpha)_n(\alpha+n)_k,
\]
the expression for Euler beta function
\[
B(\alpha,\beta_2)=\frac{\Gamma(\alpha)\Gamma(\beta_2)}{\Gamma(\alpha+\beta_2)}
\]
and the derivative formula
\begin{equation}\label{eq:2F1-derivative}
{_2F_1}'(a,b;c;x)\equiv\frac{\partial}{\partial{a}}{_2F_1}(a,b;c;x)
=\sum\limits_{k=0}^{\infty}\frac{\psi(a+k)(a)_k(b)_kx^k}{(c)_kk!}-\psi(a){_2F_1}(a,b;c;x),
\end{equation}
identity (\ref{eq:F1logarithmic}) can be rewritten as:
\begin{multline}\label{eq:F1logarithmic1}
B(\alpha,\beta_2)F_1(\alpha;\beta_1,\beta_2;\alpha+\beta_2;z_1,z_2)=
\sum\limits_{n=0}^{\infty}\frac{(\alpha)_n(\beta_2)_n}{(n!)^2}(1-z_2)^n\times
\\
\times\left\{[-\ln(1-z_2)+2\psi(1+n)-\psi(\alpha+n)-\psi(\beta_2+n)]
{_2F_1}(\alpha+n,\beta_1;\alpha;z_1)-{_2F_1}'(\alpha+n,\beta_1;\alpha;z_1)\right\}.
\end{multline}
Applying the transformation
\[
F_1(a;b_1,b_2;a+b_1+b_2;x,y)=\left(\frac{1-y}{1-x}\right)^{b1}
F_1\left(b_1+b_2;b_1,a;a+b_1+b_2;\frac{y-x}{1-x},y\right)
\]
to (\ref{eq:F1logarithmic}) and (\ref{eq:F1logarithmic1})  in view
of (\ref{eq:f-defined}) gives
\begin{multline}\label{eq:F1-mainlog}
f_{a,b_1,b_2}(x,y)=\left(\frac{1-y}{1-x}\right)^{b_1}\times
\\
\sum\limits_{k,n=0}^{\infty}\!\!\frac{(b_1)_k(b_1+b_2+k)_n(a)_n(-\ln(1-y)+2\psi(1+n)-\psi(a+n)-\psi(b_1+b_2+k+n))}{k!(n!)^2(1-y)^{-n}}
\left(\frac{y-x}{1-x}\right)^k
\\
=\left(\frac{1-y}{1-x}\right)^{b_1}
\left\{\sum\limits_{n=0}^{\infty}\frac{(a)_n(b_1+b_2)_n}{(n!)^2}\left[\ln\frac{1}{1-y}+2\psi(1+n)-\psi(a+n)-\psi(b_1+b_2+n)\right]\right.
\\
\times{_2F_1}\left(b_1+b_2+n,b_1;b_1+b_2;\frac{y-x}{1-x}\right)(1-y)^n
\\\left.
-\sum\limits_{n=0}^{\infty}\frac{(a)_n(b_1+b_2)_n}{(n!)^2}{_2F_1}'\left(b_1+b_2+n,b_1;b_1+b_2;\frac{y-x}{1-x}\right)(1-y)^n\right\}.
\end{multline}
Taking $n=0$ in the above formula and applying
\[
{_2F_1}\left(b_1+b_2,b_1;b_1+b_2;\frac{y-x}{1-x}\right)=\left(\frac{1-y}{1-x}\right)^{-b_1}
\]
we get
\begin{equation}\label{eq:b-log}
f_{a,b_1,b_2}(x,y)\!=\ln\frac{1}{1-y}
+2\psi(1)-\psi(a)-\psi(b_1+b_2)
-\left[\frac{1-y}{1-x}\right]^{b_1}\!{_2F_1}'\!\left[\!\!\begin{array}{l}b_1+b_2, b_1\\
b_1+b_2\end{array}\!\vline~\frac{y-x}{1-x}\right]+R,
\end{equation}
where it is clear from (\ref{eq:F1-mainlog}) that
\[
R=O((1-y)\ln(1-y)),
\]
which is equivalent to the second formula in (\ref{eq:R1bound}).
Formula (\ref{eq:b-log}) can be easily put into a different form
by differentiating the identity
\[
{_2F_1}(a,b;c;x)=(1-x)^{c-a-b}{_2F_1}(c-a,c-b;c;x)
\]
with respect to $a$:
\[
{_2F_1}'(a,b;c;x)=-\ln(1-x)(1-x)^{c-a-b}{_2F_1}(c-a,c-b;c;x)-(1-x)^{c-a-b}{_2F_1}'(c-a,c-b;c;x).
\]
Hence:
\[
{_2F_1}'\!\left[\!\!\begin{array}{l}b_1+b_2, b_1\\
b_1+b_2\end{array}\!\vline~\frac{y-x}{1-x}\right]=
\left(\frac{1-x}{1-y}\right)^{b_1}\ln\frac{1-x}{1-y}-\left(\frac{1-x}{1-y}\right)^{b_1}
{_2F_1}'\!\left[\!\!\begin{array}{l}0, b_2\\
b_1+b_2\end{array}\!\vline~\frac{y-x}{1-x}\right].
\]
Since
\begin{equation}\label{eq:2F1-diffat0}
F'(a,b;c;z)_{|a=0}={\sum\limits_{k=1}^{\infty}\frac{d}{da}(a)_k\frac{(b)_kz^k}{(c)_kk!}}_{|a=0}=
\sum\limits_{k=1}^{\infty}\frac{(b)_k(k-1)!}{(c)_kk!}z^k=\frac{bz}{c}{_3F_2}\left(\!\!\begin{array}{l}1,1,b+1\\2,c+1\end{array}\!\vline~z\right),
\end{equation}
we will have
\[
{_2F_1}'\!\left[\!\!\begin{array}{l}0, b_2\\
b_1+b_2\end{array}\!\vline~\frac{y-x}{1-x}\right]
=\frac{b_2(y-x)}{(b_1+b_2)(1-x)}{_3F_2}\left(\!\!\begin{array}{l}1,1,\beta_2+1\\2,\beta_1+\beta_2+1\end{array}\!\vline~\frac{y-x}{1-x}\right).
\]
In view of definition (\ref{eq:g-defined}) of $g_{a,b_1,b_2}(x,y)$
formula (\ref{eq:b-log}) transforms into (\ref{eq:g-approx}).

To estimate the remainder term we will use the ideas from
\cite{LopezMellin}. An application of the  integral representation
\cite[formula 5.8(5)]{Bat1} and a change of variable give
($u=1-x$, $v=1-y$):
\begin{multline}\label{eq:F1int}
F_1\left(a;b_1,b_2;a+b_1+b_2; 1-u,1-v\right)
\\
=u^{-b_1}v^{-b_2}
F_1\left(b_1+b_2;b_1,b_2;a+b_1+b_2;1-\frac{1}{u},1-\frac{1}{v}\right)
\\
=\frac{\Gamma(a+b_1+b_2)}{\Gamma(a)\Gamma(b_1+b_2)}
\int\limits_{0}^{\infty}\frac{t^{a-1}(1+t)^{-a}dt}{(1+ut)^{b_1}(1+vt)^{b_2}}
=\frac{\Gamma(a+b_1+b_2)}{\Gamma(a)\Gamma(b_1+b_2)}
\int\limits_{0}^{\infty}f_{a}(t)h_{b_1,b_2}(u,v;t)dt.
\end{multline}
where
\begin{equation}\label{eq:falpha-defined}
f_{a}(t)=t^{a-1}(1+t)^{-a}=\sum\limits_{k=0}^{n-1}(-1)^k\frac{(a)_k}{k!t^{k+1}}+f_{a,n}(t),~~
t\to{\infty},
\end{equation}
\begin{equation}\label{eq:f-asymp}
f_{a}(t)=O(t^{a-1}),~~t\to{0}~\Rightarrow~f\in{\cal{F}}_{1,1-a},
\end{equation}
and
\begin{equation}\label{eq:h-defined}
h_{b_1,b_2}(u,v;t)=\frac{1}{(1+tu)^{b_1}(1+tv)^{b_2}}=\sum\limits_{k=0}^{n-1}(-1)^kt^k
\sum\limits_{m=0}^{k}\frac{(b_1)_m(b_2)_{k-m}}{m!(k-m)!}u^mv^{k-m}+h_n(t),~~t\to{0};
\end{equation}
\begin{equation}\label{eq:h-asymp}
h_{b_1,b_2}(u,v;t)=O(t^{-b_1-b_2}),~t\to\infty~\Rightarrow~h\in{\cal{H}}_{0,b_1+b_2}.
\end{equation}
Spaces ${\cal F}$ and ${\cal H}$ are defined in
\cite{LopezMellin}. Basically, they mean nothing other than the
asymptotic formulas satisfied by $f$ and $h$, presented above. If
$a$, $b_1$, $b_2$ are all positive conditions I and II from
\cite{LopezMellin} are satisfied.

Representation (\ref{eq:F1int}) is not precisely a Mellin
convolution. However, if we approach the point $u=v=0$ (i.e.
$x=y=1$) along straight lines we can put $u=\gamma_1\varepsilon$,
$v=\gamma_2\varepsilon$, where $\gamma_1$ and $\gamma_2$ are
positive constants and $\varepsilon\to{0}$.  It this case
\[
h_{b_1,b_2}(u,v;t)=h_{b_1,b_2,\gamma_1,\gamma_2}(\varepsilon{t})
\]
and (\ref{eq:F1int}) takes the form of Mellin convolution. Since
every point $u$, $v$ lies on some straight line with endpoint
$(1,1)$ and all our further speculations assume sufficiently small
but fixed $u$, $v$ there are always $\gamma_1$, $\gamma_2$ and
$\varepsilon$ (of course non-unique) which are implied.  Hence the
theory from \cite{LopezMellin} can be applied.

From
\[
f_n(t)=\sum\limits_{k=n}^{\infty}(-1)^k\frac{(a)_k}{k!t^{k+1}}
=\frac{(-1)^n(a)_n{_2F_1(a+n,1;1+n;-1/t)}}{t^{n+1}n!}
=\frac{(-1)^n(a)_n}{t^{n+1}(n-1)!}\int\limits_{0}^{1}\frac{(1-s)^{n-1}}{(1+s/t)^{n+a}}ds
\]
it is obvious that $\sign(f_n)=(-1)^n$.  Similarly, from
\[
h_{b_1,b_2,n}(u,v;t)=\sum\limits_{k=n}^{\infty}(-1)^kt^k
\sum\limits_{m=0}^{k}\frac{(b_1)_m(b_2)_{k-m}}{m!(k-m)!}u^mv^{k-m}
=\sum\limits_{k=n}^{\infty}(-1)^k\frac{(b_2)_kv^kt^k}{k!}{_2F_1}(b_1,-k;1-b_2-k;u/v).
\]
it can be seen that $\sign(h_n)=(-1)^n$.  This shows that the
remainder is always positive which implies in particular
inequality (\ref{eq:fgineq}).

Now take $n=1$ and apply  \cite[Theorem~4.3]{LopezMellin} which
shows that the remainder has the form (since $a=0$, $b=1$ in terms
of \cite{LopezMellin})
\begin{multline}\label{eq:Mellinremainer}
R_{a,b_1,b_2}(u,v)=\int\limits_{0}^{\infty}f_{a,1}(t)h_{b_1,b_2,1}(u,v;t)dt
\\
=\int\limits_{0}^{\infty}
\left[\frac{t^{a-1}}{(1+t)^{a}}-\frac{1}{t}\right]\left[\frac{1}{(1+ut)^{b_1}(1+vt)^{b_2}}-1\right]dt.
\end{multline}
The bound for $R_{a,b_1,b_2}(u,v)$ is based on the following lemma
whose proof we postpone until the end of the proof of the theorem.
\begin{lemma}\label{lm:fhineq}
For all $t\in(0,\infty)$ the inequalities
\begin{equation}\label{eq:f11}
-a/t^2<f_{a,1}(t)<0,
\end{equation}
\begin{equation}\label{eq:f12}
-1/t<f_{a,1}(t)<0,
\end{equation}
\begin{equation}\label{eq:h11}
-1<h_{b_1,b_2,1}(u,v;t)<0,
\end{equation}
\begin{equation}\label{eq:h12}
-t(ub_1+vb_2)<h_{b_1,b_2,1}(u,v;t)<0
\end{equation}
hold true.
\end{lemma}
The integral in (\ref{eq:Mellinremainer}) may be decomposed as
follows
\[
R_1=\int\limits_{0}^{1}f_{a,1}(t)h_{b_1,b_2,1}(u,v;t)dt+\int\limits_{1}^{1/r}f_{a,1}(t)h_{b_1,b_2,1}(u,v;t)dt
+\int\limits_{1/r}^{\infty}f_{a,1}(t)h_{b_1,b_2,1}(u,v;t)dt,
\]
where $r$ can be any positive number (it is not needed that
$r<1$!). Set $r=ub_1+vb_2$ and use estimates (\ref{eq:f12}) and
(\ref{eq:h12}) in the first integral, (\ref{eq:f11}) and
(\ref{eq:h12}) in the second and (\ref{eq:f11}) and (\ref{eq:h11})
in the third. This gives the estimate (\ref{eq:R1bound}).
$\square$

\rem We could use Proposition~3.1 from \cite{LopezMellin}  to give
an estimate for the error term. However, in our specific situation
we are able to derive a much better bound based on
Lemma~\ref{lm:fhineq} using the method of proof of this
proposition but not it's statement.

\textbf{Proof of Lemma~\ref{lm:fhineq}.}

(a) Inequality (\ref{eq:f11}). Write $f_{a,1}(t)=g_{a}(t)/t^2$,
where
\[
g_{a}(t)=\frac{t^{a+1}}{(1+t)^{a}}-t.
\]
Then (\ref{eq:f11}) is equivalent to $-a<g_{a}(t)<0$. Clearly,
$g_{a}(0)=0$. It is an easy exercise to check that
$g_{a}(\infty)=-a$.  If we prove that $g_{a}'(t)<0$ we are done.
Differentiating and multiplying both sides by $(1+t)^{a+1}$ we see
that the required inequality takes the form
\[
(1+a)(1+t)t^{a}<(1+t)^{a+1}+a{t^{a+1}}~\Leftrightarrow~\frac{(1+t)^{a+1}}{t^{a}(1+a+t)}>1
~\Leftrightarrow~(1+x)^{a+1}>1+(1+a)x,
\]
where $x=1/t$ and the last inequality is the classical Bernoulli
inequality valid for $a>0$ and $x>-1$ \cite[formula
III(1.2)]{Mitrinovic}.

(b) Inequality (\ref{eq:f12}) is proved similarly but simpler.

(c) Inequality (\ref{eq:h11}) is obvious from the definition
(\ref{eq:h-defined}) of $h_{b_1,b_2}(u,v;t)$.

(d) To prove (\ref{eq:h12}) we again apply Bernoulli's inequality
\cite[formula III(1.2)]{Mitrinovic} in the form ($b_1,b_2>0$):
\[
(1+tu)^{-b_1}>1-b_1tu,~~~(1+tu)^{-b_2}>1-b_2tu.
\]
Multiplying these two inequalities we get the estimate
\begin{equation}\label{eq:hstrongerineq}
1-\frac{1}{(1+tu)^{b_1}(1+tv)^{b_2}}<t(ub_1+vb_2)-t^2uvb_1b_2
\end{equation}
which is even stronger than (\ref{eq:h12}). $\square$

\rem Application of (\ref{eq:hstrongerineq}) instead of
(\ref{eq:h12}) in the proof of theorem~\ref{th:fgapprox} leads to
an estimate of the remainder term $R$ which is better than
(\ref{eq:R1bound}). However, numerically it is only a very minor
improvement, so we decided to keep the simpler estimate
(\ref{eq:R1bound}) in the theorem.

\rem Representation (\ref{eq:g-defined}) also leads to the
following observation:  \emph{for general values of parameters
there exists no approximation for $f_{a,b_1,b_2}$ in the
neighbourhood of $(1,1)$ in terms of elementary functions.}
Indeed, let
\[
f_{a,b_1,b_2}(x,y)=h(x,y)+o(1)
\]
as $x,y\to{1}$ with an elementary $h(x,y)$. Then from
(\ref{eq:g-approx}):
\[
g_{a,b_1,b_2}(x,y)-h(x,y)=o(1)~\Rightarrow
\]\[
h(x,y)+\ln(1-x)=
2\psi(1)-\psi(a)-\psi(b_1+b_2)+\frac{b_2(y-x)}{(b_1+b_2)(1-x)}
{_3F_2}\left(\!\!\begin{array}{l}1,1,b_2+1\\2,b_1+b_2+1\end{array}\!\vline~\frac{y-x}{1-x}\right)+\varepsilon(x,y),
\]
and $\varepsilon(x,y)\to{0}$ as $x,y\to{1}$. Let  $x,y\to{1}$
along a straight line going through $(1,1)$, so that
\[
(1-y)/(1-x)=\gamma=\mathrm{const}.
\]
Then, due to
\[
(y-x)/(1-x)=1-\gamma,
\]
we have for the elementary $h_1(x,y)\equiv h(x,y)+\ln(1-x)$:
\[
h_1(x,y)= d(\gamma)+\varepsilon(x,y),
\]
where
\[
d(\gamma)=2\psi(1)-\psi(a)-\psi(b_1+b_2)+\frac{b_2(1-\gamma)}{(b_1+b_2)}
{_3F_2}\left(\!\!\begin{array}{l}1,1,b_2+1\\2,b_1+b_2+1\end{array}\!\vline~1-\gamma\right).
\]
Since $y=1-(1-x)\gamma$,  we can write the above as
\[
\tilde{h}_1(x,\gamma)=d(\gamma)+\tilde{\varepsilon}(x,\gamma),
\]
where $\gamma\in(0,\infty)$ is arbitrary, but fixed.  For $x=1$
this gives $\tilde{h}_1(1,\gamma)=d(\gamma)$ for all
$\gamma\in(0,\infty)$.  Hence, a restriction of an elementary
function $\tilde{h}_1$ gives ${_3F_2}$ for all values of its
argument in the range $(-\infty,1)$, which is impossible, and so
$h(x,y)$ cannot be an elementary function.

\rem Expansion \cite[formula (53)]{Fer-Lopez} can ce cast into the
form
\begin{multline}\label{eq:main-lopez}
\frac{\Gamma(b_1+b_2)\Gamma(a)}{\Gamma(a+b_1+b_2)}F_1\left(a;b_1,b_2;a+b_1+b_2;
1-\frac{\gamma_1}{z},1-\frac{\gamma_2}{z}\right)=
\\
\sum\limits_{k=0}^{n-1}\left[\frac{D_k(a,b_1,b_2;\gamma_1,\gamma_2)}{z^k}
+\log(z)\frac{E_k(a,b_1,b_2;x,y)}{z^k}\right]+R_n(a,b_1,b_2,\gamma_1,\gamma_2;z),
\end{multline}
Substituting $x=1-\gamma_1/z$, $y=1-\gamma_2/z$ into
(\ref{eq:F1-mainlog}) we see that both (\ref{eq:main-lopez}) and
(\ref{eq:F1-mainlog}) are asymptotic expansions for $|z|\to\infty$
in the same asymptotic sequences $z^{-k}$, $z^{-k}\log(z)$ and so
their coefficients are the same.  Hence, (\ref{eq:F1-mainlog}) can
be viewed as a simpler form of \cite[formula (53)]{Fer-Lopez}. The
appearance of the coefficients $D_k$ and $E_k$ is very different
from that of the coefficients of (\ref{eq:F1-mainlog}) and direct
reduction is non-trivial. For instance, the first term of
\cite[formula (53)]{Fer-Lopez} reads (after some simple
manipulations) ($F={_2F_1}$, $M=(1-y)/(1-x)$):
\begin{multline}
B(a,b_1+b_2)F_1(a;b_1,b_2;a+b_1+b_2;x,y)=\psi(1)-\psi(a)\\
+\frac{-\ln(1-v)+\ln(M)+\psi(1)-\psi(b_1+b_2)}{b_1+b_2}\left(Mb_2
F\left[\!\!\begin{array}{l}1,b_2+1\\b_1+b_2+1\end{array}\!\vline~1-M\right]+
b_1F\left[\!\!\begin{array}{l}1,b_2\\b_1+b_2+1\end{array}\!\vline~1-M\right]\right)\\
+\frac{1}{b_1+b_2}\left(Mb_2
F'\left[\!\!\begin{array}{l}1,b_2+1\\b_1+b_2+1\end{array}\!\vline~1-M\right]+
b_1F'\left[\!\!\begin{array}{l}1,b_2\\b_1+b_2+1\end{array}\!\vline~1-M\right]\right)+R_1.
\end{multline}
Now using the relation \cite[formula 2.8(36)]{Bat1}
\begin{equation}\label{eq:2F1rec}
(c-a-b)F(a,b;c;z)-(c-a)F(a-1,b;c;z)+b(1-z)F(a,b+1;c;z)=0
\end{equation}
we immediately get
\[
Mb_2F\left[\!\!\begin{array}{l}1,b_2+1\\b_1+b_2+1\end{array}\!\vline~1-M\right]+
b_1F\left[\!\!\begin{array}{l}1,b_2\\b_1+b_2+1\end{array}\!\vline~1-M\right]=b_1+b_2.
\]
Differentiating (\ref{eq:2F1rec}) with respect to $a$ and putting
$a=0$ we get:
\[
(c-b-1)F'(1,b;c;z)+b(1-z)F'(1,b+1;c;z)=F(1,b;c;z)+(c-1)F'(0,b;c;z)-1.
\]
Using (\ref{eq:2F1-diffat0}) we see
\[
Mb_2F'\left[\!\!\begin{array}{l}1,b_2+1\\b_1+b_2+1\end{array}\!\vline~1-M\right]+
b_1F'\left[\!\!\begin{array}{l}1,b_2\\b_1+b_2+1\end{array}\!\vline~1-M\right]
\]\[
=F\left[\!\!\begin{array}{l}1,b_2\\b_1+b_2+1\end{array}\!\vline~1-M\right]
+\frac{b_2(b_1+b_2)(1-M)}{(b_1+b_2+1)}
{_3F_2}\left[\!\!\begin{array}{l}1,1,b_2+1\\2,b_1+b_2+2\end{array}\!\vline~1-M\right]-1
\]
and
\begin{multline}\label{eq:firstorder}
B(a,b_1+b_2)F_1(a;b_1,b_2;a+b_1+b_2;x,y)=
\ln\frac{1}{1-x}+2\psi(1)-\psi(a)-\psi(b_1+b_2)
\\+\frac{1}{b_1+b_2}F\left[\!\!\begin{array}{l}1,b_2\\b_1+b_2+1\end{array}\!\vline~1-M\right]
+\frac{b_2(1-M)}{(b_1+b_2+1)}
{_3F_2}\left[\!\!\begin{array}{l}1,1,b_2+1\\2,b_1+b_2+2\end{array}\!\vline~1-M\right]-\frac{1}{b_1+b_2}+R_1.
\end{multline}
Finally, (\ref{eq:firstorder}) is reduced to (\ref{eq:g-defined})
with
the help of the following formula found at\\
http://functions.wolfram.com/07.27.17.0029.01:
\[
{_3F_2}(a,b,c;a+1,e;z)=\frac{1}{a-e+1}\left[a{_2F_1}(b,c;e;z)-(e-1){_3F_2}(a,b,c;a+1,e-1;z)\right].
\]
Recalling that $M=(1-y)/(1-x)$ we get (\ref{eq:g-defined}).  The
direct reduction for further terms is even more complicated.

\begin{theo} The following properties are true:
\begin{enumerate}
\item The function $g$ is permutation symmetric:
\begin{equation}\label{eq:gsymmetry}
g_{a,b_1,b_2}(x,y)=g_{a,b_2,b_1}(y,x).
\end{equation}

\item For $y=1$ \emph{(}and $x=1$ due to
\emph{(\ref{eq:gsymmetry})}\emph{)} the function
$g_{a,b_1,b_2}(x,y)$ reduces to the Ramanujan's approximation:
\begin{subequations}\label{eq:gy1}
\begin{align}
g_{a,b_1,b_2}(x,1)=\ln\frac{1}{1-x}+2\psi(1)-\psi(a)-\psi(b_1),\\
g_{a,b_1,b_2}(1,y)=\ln\frac{1}{1-y}+2\psi(1)-\psi(a)-\psi(b_2).
\end{align}
\end{subequations}

\item For $x=y$ the function $g_{a,b_1,b_2}(x,y)$ becomes the
Ramanujan's approximation:
\begin{equation}\label{eq:g-diagonal}
g_{a,b_1,b_2}(x,x)=\ln\frac{1}{1-x}
+2\psi(1)-\psi(a)-\psi(b_1+b_2).
\end{equation}

\item For the values of parameters $a=b_1=b_2=1/2$ we have
\begin{equation}\label{eq:gCG1}
g_{1/2,1/2,1/2}=\ln\frac{4}{\sqrt{1-\lambda^2}+\sqrt{1-k^2\lambda^2}},
\end{equation}
which is the approximation of Carlson-Gustafson.
\end{enumerate}
\end{theo}

\textbf{Proof.} To prove the first statement we need the following
elementary lemma:
\begin{lemma}
For $b\ne{1}$ the following relation holds true:
\begin{equation}\label{eq:3F2relation}
{_3F_2}\left(\!\!\begin{array}{l}1,b,c\\2,e\end{array}\!\vline~\frac{z}{z-1}\right)=
\frac{(1-z)^b(c-e)}{c-1}{_3F_2}\left(\!\!\begin{array}{l}1,b,e-c+1\\2,e\end{array}\!\vline~z\right)+
\frac{(e-1)(1-z)(1-(1-z)^{b-1})}{(c-1)(b-1)z}.
\end{equation}
For $b=1$ it reduces to
\begin{equation}\label{eq:3F2logrelation}
{_3F_2}\left(\!\!\begin{array}{l}1,1,c\\2,e\end{array}\!\vline~\frac{z}{z-1}\right)=
\frac{(z-1)(e-c)}{c-1}{_3F_2}\left(\!\!\begin{array}{l}1,1,e-c+1\\2,e\end{array}\!\vline~z\right)+
\frac{(e-1)(1-z)}{(c-1)z}\ln\frac{1}{1-z}.
\end{equation}
\end{lemma}

\textbf{Proof.} The proof is based on the following easily
verifiable relation (which can be also found at
http://functions.wolfram.com/07.27.03.0120.01):
\begin{equation}\label{eq:3F2to2F1}
{_3F_2}\left(\!\!\begin{array}{l}1,b,c\\2,e\end{array}\!\vline~z\right)=\frac{e-1}{(b-1)(c-1)z}
\left[{_2F_1}\left(\!\!\begin{array}{l}b-1,c-1\\e-1\end{array}\!\vline~z\right)-1\right].
\end{equation}
To prove (\ref{eq:3F2relation}) write this relation for $z/(z-1)$
in place of $z$, apply
\[
{_2F_1}\left(\!\!\begin{array}{l}b-1,c-1\\e-1\end{array}\!\vline~\frac{z}{z-1}\right)
=(1-z)^{b-1}{_2F_1}\left(\!\!\begin{array}{l}b-1,e-c\\e-1\end{array}\!\vline~z\right)
\]
and substitute ${_2F_1}$ from the right-hand side by ${_2F_1}$
expressed from (\ref{eq:3F2to2F1}).  To prove
(\ref{eq:3F2logrelation}) let $b$ tend to $1$ and apply the
L'Hopital rule.~$\square$

Combining (\ref{eq:3F2logrelation}) with the definition
(\ref{eq:g-defined}) of $g_{a,b_1,b_2}(x,y)$  we immediately
obtain (\ref{eq:gsymmetry}).

Next we check the behavior of  the function $g_{a,b_1,b_2}(x,y)$
on the sides of the square $|x|<1$, $|y|<1$. Writing
(\ref{eq:3F2to2F1}) for $z=1$ and using the Gauss formula for
${_2F_1(1)}$ we get
\begin{multline*}
{_3F_2}\left(\!\!\begin{array}{l}1,b,c\\2,e\end{array}\!\vline~1\right)=\frac{e-1}{(b-1)(c-1)}
\left[{_2F_1}\left(\!\!\begin{array}{l}b-1,c-1\\e-1\end{array}\!\vline~1\right)-1\right]
\\
=\frac{e-1}{(b-1)(c-1)}\left[\frac{\Gamma(e-1)\Gamma(e-b-c+1)}{\Gamma(e-b)\Gamma(e-c)}-1\right].
\end{multline*}
 Now let $b\to{1}$ and use the L'Hopital rule:
\begin{multline*}
{_3F_2}\left(\!\!\begin{array}{l}1,1,c\\2,e\end{array}\!\vline~1\right)
=\left.\frac{(e-1)\Gamma(e-1)}{(c-1)\Gamma(e-c)}\frac{d}{db}\frac{\Gamma(e-b-c+1)}{\Gamma(e-b)}\right|_{b=1}
\\
=\frac{\Gamma(e)}{(c-1)\Gamma(e-c)}\frac{-\Gamma(e-c)\psi(e-c)\Gamma(e-1)+\Gamma(e-1)\psi(e-1)\Gamma(e-c)}{[\Gamma(e-1)]^2}
=\frac{(e-1)}{(c-1)}(\psi(e-1)-\psi(e-c)).
\end{multline*}
Substituting $e=b_1+b_2+1$, $c=b_2+1$ gives (\ref{eq:gy1}).

Identity (\ref{eq:g-diagonal}) is obvious from the definition
(\ref{eq:g-defined}) of $g_{a,b_1,b_2}(x,y)$.

Finally, formula (\ref{eq:gCG1}) follows from the reduction
formula
\[
{_3F_2}\left(\!\!\begin{array}{l}1,1,3/2\\2,2\end{array}\!\vline~z\right)
=-\frac{4}{z}\ln\left(\frac{1}{2}+\frac{\sqrt{1-z}}{2}\right).
\]
This completes the proof of the theorem. $\square$
\begin{corolt}\label{cl:rough} For $x,y\to{1}$
\begin{equation}\label{eq:f-roughapprox}
f_{a,b_1,b_2}(x,y)=\ln\frac{1}{1-xy}+O(1).
\end{equation}
\end{corolt}

\textbf{Proof.}  Assume first that $x$ and $y$ approach $(1,1)$ in
a way such $(1-y)/(1-x)$ stays bounded.  We have
\[
\ln\frac{1}{1-xy}=\ln\frac{1}{1-x+x-xy}=\ln\frac{1}{(1-x)\left(1+x\frac{1-y}{1-x}\right)}=\ln\frac{1}{1-x}+\ln\frac{1}{1+x\frac{1-y}{1-x}}.
\]
Hence,
\begin{multline*}
\ln\frac{1}{1-xy}-g_{a,b_1,b_2}(x,y)
\\
=\ln\frac{1}{1+x\frac{1-y}{1-x}}-\gamma(a,b_1+b_2)
-\frac{b_2(y-x)}{(b_1+b_2)(1-x)}
{_3F_2}\left(\!\!\begin{array}{l}1,1,b_2+1\\2,b_1+b_2+1\end{array}\!\vline~\frac{y-x}{1-x}\right)=O(1).
\end{multline*}
If $(1-y)/(1-x)$ is unbounded, than exchange the roles of $x$ and
$y$ and use (\ref{eq:gsymmetry}).~$\square$

Finally, we remark that the authors of \cite{ABRVV,AVV} consider
monotonicity and ranges of the functions
\[
\frac{1-{_2F_1}(a,b;a+b;x)}{\ln(1-x)},
~~\frac{x{_2F_1}(a,b;a+b;x)}{\ln(1/(1-x))}
\]
and
\[
B(a,b){_2F_1}(a,b;a+b;x)+\ln(1-x)
\]
for $x\in(0,1)$.  Our Corollary~\ref{cl:rough}  shows that similar
problems can be considered for the combinations
\[
\frac{1-F_1(\alpha;\beta_1,\beta_2;\alpha+\beta_1+\beta_2;x,y)}{\ln(1-xy)},
\]
and
\[
f_{\alpha,\beta_1,\beta_2}(x,y)-\ln\frac{1}{1-xy}
\]
for $x,y\in(0,1)$.

\paragraph{3. Acknowledgments.} The author is thankful to Professor
J.L.\,L\'{o}pez of Universidad P\'{u}blica de Navarra and
Professor Matti Vuorinen of University of Helsinki for a series
for useful discussions.  This research has been supported by INTAS
(grant no.05-109-4968),  the Russian Basic Research Fund (grant
no. 05-01-00099) and Far Eastern Branch of the Russian Academy of
Sciences (grant no. 06-III-B-01-020).


\begin{thebibliography}{9}
\bibitem{ABRVV} G.D.\,Anderson, R.W.\,Barnard,  K.C.\,Richards,
M.K.\,Vamanamurthy and M.\,Vuorinen, Inequalities for
zero-balanced hypergeometric functions. \emph{Trans. Amer. Math.
Soc.} 347 (1995), no. 5, 1713--1723.
\bibitem{AVV}G.D.\,Anderson, M.K.\,Vamanamurthy and M.\,Vuorinen,
\emph{Conformal invariants, inequalities, and quasiconformal
maps}. Canadian Mathematical Society Series of Monographs and
Advanced Texts. A Wiley-Interscience Publication. John Wiley and
Sons, Inc., New York, 1997.
\bibitem{Berndt0}B.C.\,Berndt, Chapter 11 of Ramanujan's second notebook, \emph{Bull. London Math.
Soc.} 15(1983), 273-320.
\bibitem{Berndt1}B.C.\,Berndt, \emph{Ramanujan's notebooks. Part
I}, Springer-Verlag, New York, 1985.
\bibitem{Berndt2} B.C.\,Berndt, \emph{Ramanujan's Notebooks,
Part II}, Springer-Verlag, New York, 1989.
\bibitem{Carlson} B.C.\,Carlson, Some series and bounds for incomplete elliptic integrals, \emph{J.
 Math. Phys.} {\bf 40}(1961), 125-134.
\bibitem{Buhring}W. B\"{u}hring, Generalized hypergeometric functions at unit argument, \emph{Proc. Amer.
Math. Soc.} 114(1992), 145-153.
\bibitem{Srivastava}W.\,B\"{u}hring and H.M.\,Srivastava, Analytic
continuation of the generalized hypergeometric series near unit
argument with emphasis on the zero-balanced series.
\emph{Approximation theory and applications}, 17--35, Hadronic
Press, Palm Harbor, FL, 1998.
\bibitem{CG1} B.C.\,Carlson and J.L.\,Gustafson Asymptotic
expansion of the first elliptic integral, \emph{SIAM J.
Math.Anal.}, vol.\textbf{16} (1985), no.5, 1072-1092.
\bibitem{Bat1} Erd\'{e}lyi~A., Magnus~W., Oberhettinger~F. and
Tricomi~F.G., 1953, \emph{Higher transcendental functions, Vol. 1}
(New York: McGraw-Hill Book Company, Inc.).
\bibitem{Fer-Lopez} C.\,Ferreira and J.\,L.\,L\'{o}pez, Asymptotic
Expansions of the Appell's Function $F_1$, \emph{Q. Appl. Math.}
62, No.2, 235-257 (2004).
\bibitem{Ismail} M.E.H.\,Ismail and J.\,Pitman, Algebraic Evaluations
of Some Euler Integrals, Duplication Formulae for Appell's
Hypergeometric Function F1, and Brownian Variations, \emph{Canad.
J. Math.} Vol. 52 (5), 2000 pp. 961-981.
\bibitem{KarpSitnik}  D.\,Karp and S.M.\,Sitnik, Asymptotic approximations
for the first incomplete elliptic integral near logarithmic
singularity, \emph{J. of Comp. and Appl. Math.}, in press.
Available on-line via http://dx.doi.org/10.1016/j.cam.2006.04.053
\bibitem{LopezMellin} J.L.\,L\'{o}pez, Asymptotic expansion of
Mellin convolution integral, JCAM, submitted.
\bibitem{Marichev} O.I.\,Marichev and S.L.\,Kalla, Behaviour of
hypergeometric function ${_pF_{p-1}(z)}$ in the vicinity of unity,
\emph{Rev. T\'{e}cn. Fac. Ingr. Univ. Zulia}, 7(1984), 1-8.
\bibitem{Mitrinovic}D.S.\,Mitrinovic, J.E.\,Pecaric, A.M.\,Fink, \emph{Classical and new
inequalities in Analysis.} Kluwer Academic Publishers, 1993.
\bibitem{Norlund} N.\,N$\o$rlund, Hypergeometric functions, \emph{Acta Math.}, 94(1955),
289-349.
\bibitem{Saigo}M.\,Saigo and H.M.\,Srivastava,
The behavior of the zero-balanced hypergeometric series ${}\sb
pF\sb {p-1}$ near the boundary of its convergence region.
\emph{Proc. Amer. Math. Soc.} 110 (1990), no. 1, 71-76.

\end{thebibliography}
\end{document}